%% file: Isometric_Actions_and_Finite_approximations_ETDS.tex
\documentclass[12pt]{article}

\input{preamble}

\begin{document}
\title{Isometric Actions and Finite Approximations}
\author{SJ Pilgrim}
\affil{University of Hawai'i at M\={a}noa}
\maketitle

\paragraph{Abstract:}  
We show that every isometric action on a Cantor set is conjugate to an inverse limit of actions on finite sets; and that every isometric action by a finitely generated amenable group is residually finite.  


\paragraph{MSC classes:} 37B99 (primary), 47L65 (secondary), 20F65 (secondary)

\section{Introduction}

Finite (dimensional) approximation properties are ubiquitous in the study of operator algebras.  This in part motivated a formulation of residual finiteness for group actions, wherein a topological action is approximated by actions on finite sets.  This note documents part of the relationship between actions which fix a metric and actions which admit such finite approximations.  In particular, it shows that every isometric action on a Cantor set is conjugate to an inverse limit of actions on finite sets; and that every isometric action by a free product of amenable groups is residually finite.  

These results give evidence for the possibility of obtaining new information about operator algebras by providing approximations in dynamics, as residual finiteness relates to finite dimensional approximation properties (namely MF-ness and quasidiagonality) when one passes to the crossed product $C^*$-algebra.  Future consideration of finite approximation properties of group actions may also yield interesting results within the study of dynamics more properly, as the former of the theorems mentioned above hints at.  This work was partly supported by the NSF (DMS 1901522).  



\definition{An action $\Gamma\curvearrowright X$ on a compact space $X$ by homeomorphisms is \textit{residually finite} if for every finite $F\subset \Gamma$ and every neighborhood $\epsilon$ of the diagonal in $X\times X$, there exist a finite set $E$, an action $\Gamma\curvearrowright E$ and a map $\zeta: E\to X$ such that $(\zeta(f\cdot e), f\cdot \zeta(e))\in \epsilon$ for all $f\in F$ and $\zeta(E)$ is $\epsilon$-dense in $X$\footnote{This means for any $x\in X$ that $\epsilon(x) = \{y : (x, y)\in \epsilon\}$ contains some $\zeta(e)$.}.  We will refer to $\Gamma\acts E$ together with $\zeta$ as a \textit{residue} of $\Gamma\acts X$ and say it is $(\epsilon, F)$-approximating.  We call a net of $(\epsilon_\alpha, F_\alpha)$-approximating residues $(\Gamma\acts E_\alpha)_{\alpha\in A}$ with $\bigcup_\alpha F_\alpha = \Gamma$ an increasing union and $\{\epsilon_\alpha\}$ a decreasing neighborhood basis of the diagonal a \textit{filtration}.  Note that it is enough to only assume the $\epsilon$-density condition locally: if for all $x\in X$ and $\epsilon>0$ there is a residue with some point within $\epsilon$ of $x$, we can apply this fact inside each ball in a finite cover of $X$ by $\epsilon$-balls to get an action satisfying the original definition.  }  

This definition was first introduced in \cite{kerr2011residually} and we recall here a few basic facts proved in that paper (immediately after definition 2.1).  First, if $X$ has no isolated points, a perturbation argument shows we can take $\zeta$ to be an inclusion map. Next, a group admits a free, residually finite action if and only if it is residually finite.  Finally, every residually finite action fixes a measure.

\normalfont

A finitely generated group acting isometrically on a compact space is residually finite.  This is a consequence of the Arzel\`{a}-Ascoli theorem, which shows that such a group embeds in a compact group; together with the Peter-Weyl theorem (which implies such a group is residually linear) and Mal'cev's theorem.  One might then ask if all isometric actions are residually finite.  Before discussing this question, one remark is in order.  

\remark{If $\Gamma\acts X$ is equicontinuous in the sense that the image of $\Gamma$ in $\Aut(X)$ is equicontinuous, then $\Gamma\acts X$ fixes a compatible metric on $X$: the image of $\Gamma$ in $\Aut(X)$ is precompact by Arzel\`{a}-Ascoli and so we can average the original metric on $X$ over the closure of the image of $\Gamma$.  Since residual finiteness of group actions and other topological properties do not depend on the metric on $X$, the theorems of proved here in fact work for equicontinuous actions.  } \label{equicontinuous remark}

\section{Actions on Cantor sets}

\normalfont
Isometric actions on Cantor sets turn out to have a finite approximation property stronger than residual finiteness (and similar to a different definition of the same term in \cite{automata}), which  leads to a structural result.

\definition{A \textit{Cantor system} is a topological dynamical system $\Gamma\acts X$ where $X$ is homeomorphic to the Cantor set.  Such a system is \textit{isometric} if it fixes a compatible metric on $X$.  }

\definition{We will call a sequence of $\Gamma$-actions with surjective, equivariant connecting maps an $\textit{inverse sequence}$ of $\Gamma$-actions.  Such a sequence has an inverse limit constructed in the usual way.  An \textit{odometer} is an action $\Gamma\acts \bf{E}$ which is an inverse limit of actions on finite sets.  We call the action on a given finite set a \textit{level} of the odometer.  Notice that an odometer is minimal iff the action on each level is transitive.  }

Observe that odometers obey a finite approximation property stronger than residual finiteness in that the estimate is not relative to a finite $F\subset \Gamma$; and also that every odometer fixes a metric (we will show the converse of this in \ref{cantor systems theorem} ).  A similar definition of `odometer' arises in \cite{palindromicsubshifts}, as does a definition of `residually finite' for actions on Cantor sets; and a structural result equivalent to \ref{cantor systems theorem} is shown in \cite[6.4]{automata} for Cantor systems which are residually finite in the sense of \cite{palindromicsubshifts}.  The proofs here will show incidentally that isometric Cantor systems satisfy this version of residual finiteness as well.  

The anonymous referee pointed out the following result is known to some experts, but it does not appear to be in the literature: for example, in the minimal case, the result can be found in the PhD thesis of Jessica Dyer.

\normalfont


\theorem{Suppose $\Gamma\acts X$ is a Cantor system (or more generally a topological action on a compact, metrizable, totally disconnected space).  Then $\Gamma\acts X$ is conjugate (isomorphic) to an inverse limit of $\Gamma$-actions on finite sets (hence an odometer) if and only if it is equicontinuous.  If the action is minimal, each level of the odometer is a transitive action.  In particular,  isometric Cantor systems are residually finite (acording to the definition above \textit{and} that of \cite{palindromicsubshifts}).  } \label{cantor systems theorem}

\begin{proof}

By \ref{equicontinuous remark}, we can assume $\Gamma\acts X$ is isometric.  

For $R>0$, define the $R$-connected components of $X$ to be the maximal sets which cannot be written as a disjoint union of non-empty open sets which are distance $\geq R$ from each other.  Since each such component is at least distance $R$ from its complement in $X$, each is a clopen set.  Notice that belonging to the same $R$-connected component of $X$ is an equivalence relation on $X$ and that the action $\Gamma\acts X$ permutes these equivalence classes.  


Consider the decomposition of $X$ into its $R$-connected components with $R<\epsilon$ and the action of $\Gamma$ on these equivalence classes induced from $\Gamma\acts X$.  



The decomposition described above is a disjoint cover of $X$.  If we repeat all of this again with $\epsilon'<\epsilon$, the sets in the new cover will be contained inside the sets in the original.  Building a system of covers this way and passing to the nerves (which are really just finite sets), we get that $\Gamma\acts X$ is conjugate (by way of the homeomorphism from $X$ to the inverse limit of nerves) to an inverse limit of actions on the nerves where the equivariant connecting maps come from inclusions.  Finally, if the original action was minimal, so is the odometer it is conjugate to, and an odometer is minimal iff the action on each level is transitive.  \end{proof}


\normalfont

\normalfont

\remark{This gives a fairly precise vindication of the claim in \cite[2.4]{palindromicsubshifts} that residual finiteness (as defined there) is in some sense opposite to expansiveness.  }


\normalfont

Next we discuss some consequences of the above result.  The first is an immediate consequence of combining with \cite[3.4 and 3.5]{kerr2011residually}.  

\corollary{Suppose $\Gamma\acts X$ is an equicontinuous Cantor system.  Then $C(X)\rtimes\Gamma$ is MF if $C_r^*(\Gamma)$ is MF (where the latter denotes the reduced $C^*$-algebra of $\Gamma$) and is quasi-diagonal if $\Gamma$ is amenable (equivalently if $C_r^*(\Gamma)$ is quasi-diagonal).  }\label{cantor systems corollary 1}\qed

\normalfont

We can also obtain a structural result about the crossed products arising from isometric Cantor systems.  We leave the standard proofs of the following facts to the reader.  

\lemma{Suppose $\Gamma\acts X_n$ is an inverse sequence of actions (that is, there are equivariant, surjective maps $X_{n+1}\to X_n$).  Then the crossed products $C(X_n)\rtimes \Gamma$ form a direct sequence (that is, there are injective homomorphisms $C(X_n)\rtimes \Gamma\to C(X_{n+1})\rtimes \Gamma$) and in fact $\displaystyle{\lim_\rightarrow C(X_n)\rtimes \Gamma} \cong C(\displaystyle{\lim_\leftarrow X_n})\rtimes \Gamma$.  } \qed

\lemma{Suppose $\Gamma\acts E$ is an action on a finite set.  Then the crossed product $\C^E\rtimes \Gamma$ has the form $\bigoplus_{\Gamma e\in \Gamma\text{\textbackslash} E} M_{|\Gamma e|}(C^*(\Gamma_e))$ where $\Gamma e$ is the orbit of $e$ and $\Gamma_e$ is the stabilizer group of $e$.  } \qed

\normalfont
From these standard facts and \ref{cantor systems theorem}, we immediately deduce the following corollary (which also gives a more elementary proof of \ref{cantor systems corollary 1}).  

\corollary{Suppose $\Gamma\acts X$ is an equicontinuous Cantor system.  Then there is an inverse system of finite sets $(E_n, \zeta_n)$ on which $\Gamma$ acts such that $C(X)\rtimes \Gamma$ is isomorphic to a direct limit $\displaystyle{\lim_{\rightarrow n}}\bigoplus_{\Gamma e\in \Gamma \text{\textbackslash} E_n} M_{|\Gamma e|}(C^*(\Gamma_{e, n}))$ where each $\Gamma_{e, n}$ is a finite-index subgroup of $\Gamma$ and $\Gamma_{e, n} > \Gamma_{e', n+1}$ for any $e'\in \zeta_n^{-1}(\{e\})$.  }\qed

\normalfont
Finally, we can also use the original structure theorem to prove a stability property for certain topological actions wherein isometric `almost actions' can be perturbed slightly to give actions.  The property shown here is spiritually similar to weak $C^*$-stability of a group $\Gamma$ as in \cite{eilers} which is equivalent to weak semiprojectivity of the full group $C^*$-algebra of $\Gamma$, a property sometimes relevant in the study of finite-dimensional approximations.  

\corollary{Suppose $X$ is a Cantor set with metric $d$ and $\Gamma$ is finitely-presented.  Suppose further that there are set maps $\mu_i: \Gamma\to \text{Isom}(X)$ such that for all $\epsilon>0$ and $S\subset \Gamma$ finite, there is $j$ such that $d(\mu_i(\gamma\delta)(x), \mu_i(\gamma)\circ\mu_i(\delta)(x))<\epsilon$ for all $x\in X$, $\gamma, \delta\in S$, and $i\geq j$ .  Then there are homomorphisms $\tilde{\mu}_i: \Gamma\to \text{Aut}(X)$ such that for all $\epsilon>0$ and $S\subset \Gamma$ there is $j$ such that $d(\tilde{\lambda}_i(\gamma)(x), \lambda_i(\gamma)(x))<\epsilon$ for $x\in X$, $\gamma\in S$, and $i\geq j$.  }

\begin{proof}

Suppose $\mu_i: \Gamma\to \text{Isom}(X)$ is as above.  Then each $\mu_i$ gives an isometric action $F\acts X$ by the free group on $S$ where $S$ is a finite generating set of $\Gamma$.  But this action is conjugate to an inverse limit of $\Gamma$-actions on finite sets, and so gives actions $\alpha_{N, i}$ of $F$ on $E_n$ for each finite set $E_n$.  For any fixed $N$, we can then find $i_N$ sufficiently large so that $\alpha_{N, i}$ descends to an honest action $\Gamma\acts E_N$.  

Let $j_N: \displaystyle{\lim_{n\leftarrow} E_n}\to E_N$ be the usual projection.  Then the pullbacks $j_N^{-1}(\{e\})$ for $e\in E_N$ are all homeomorphic to each other.  So we can define an action $\Gamma\acts X\cong \displaystyle{\lim_{n\leftarrow} E_n}$ by homeomorphisms (just by fixing an identification $\displaystyle{\lim_\leftarrow}E_n \cong E_N\times \mathcal{C}$ and acting by $\alpha_N\times \text{id}_\mathcal{C}$, where $\mathcal{C}$ is a Cantor set).  Call this action $\tilde{\mu}_i$.  Then for $i\geq i_N$, $\tilde{\mu}_i$ will be will be within $\epsilon_N$ of $\mu_i$ for some $\epsilon_N\to 0$ as $N\to \infty$.  \end{proof}

\section{Isometric actions by amenable groups}
\normalfont

As explained earlier, the representation theory of compact groups shows that a finitely generated subgroup of a compact group is residually finite.  One might then ask if an analogue of this result for residually finite actions holds, i.e. are translation actions by such subgroups residually finite?  Indeed we will show this when $\Gamma$ is amenable\footnote{Residual finiteness of actions by free groups preserving measures of full support was previously shown in \cite{kerr2011residually}}, and the representation theory of compact groups plays a major role.  At the end of this section, we will see this is enough to show that any isometric action by an amenable group is residually finite.  

We first need some lemmas about isometries of flat tori, that is, those which arise as quotients $\R^n/L$ where $L$ is a lattice.  

\lemma{Let $T^n$ be a torus with a flat metric.  Then the isometry group of $T^n$ is isomorphic to $F\ltimes T^n$ for a finite group $F$ and if $E^n$ is the image of $\frac{1}{N}L$ in $\mathbb{R}^n/L$ for some $N$ where $L$ is the lattice used to construct $T^n$, then any isometry coming from $F$ takes $E^n$ into itself.}

\begin{proof}

Let $g\in \text{Isom}(T^n)$.  Let $h\in Isom(T^n)$ be given by $h(x) = g(x) - g(0)$, so $h(0) = 0$.  Since $h$ is basepoint-preserving, it lifts to an isometry $\tilde{h}$ of $\R^n$ that fixes $0\in \R^n$, i.e. an element of $O_n(\R)$.  As $\tilde{h}$ lifts a map of $T^n$, $\tilde{h}(L) = L$, so $\tilde{h}\in O_n(\R)\cap \text{Aut}(L)$ which is a finite group $F$.  Hence $h$ acts on $T^n$ via an element of $F$.  

Since we can decompose isometries of $T^n$ as an element of $F$ mulitplied by a translation (as we did with $g$), there is a surjective map $F\times T^n\to \text{Isom}(T^n)$ given by $(f, z)\mapsto (x\mapsto f\cdot x + z)$.  It is now straightforward to check that this map is injective, continuous, and a homomorphism for the semidirect product group structure on $F\ltimes T^n$.  \end{proof}

\lemma{Let $T^n$ be as in the previous lemma.  For all $\epsilon>0$ and $S$ a finite collection of the isometries discussed above, there is a finite, $\epsilon$-dense subset of $T^n$ admitting an $(\epsilon, S)$-approximating action for $\langle S\rangle \acts T^n$.  }

\begin{proof}

Let $\delta >0$.  We can assume $S$ has the form $F\cup Z$ where $F$ is the finite part of the semidirect product forming $\text{Isom}(T^n)$ as above and $Z$ is a finite collection of rotations each in a single coordinate.  Find $N$ so that each irrational rotation in $Z$ can be approximated to within $\delta$ by rational rotations with denominator $N$, and so that $||b_i/N||<\delta$ for all elements $b_i$ of some fixed basis of the lattice used to construct $T^n$.  Also make sure $N$ contains as a factor the denominator of any rational rotations in $Z$ (so that our approximation can be exact on any such rotations).  Then each of these rational rotations acts on a set $E^n$ from the previous lemma, giving a $\delta$-approximating action (for the generators in $Z$) on an $n\sqrt{n}\delta$-dense subset of $T^n$ which in fact coincides with the usual action on $E^n$ by some rational rotations inside $T^n$.  

The previous lemma now allows us to create an approximating action for the rest of $\langle S\rangle$ which is implemented through a translation action of $F\ltimes(\oplus_{k=1}^n \Z/N) < F\ltimes T^n$, as it comes from $F$ and a bunch of rational rotations with denominator $N$.  This action therefore respects the relations on the semi-direct product.  Letting $\delta$ be sufficiently small compared to $\epsilon$ finishes the proof.  \end{proof}

\lemma{Suppose $\Gamma<G$ is a finitely-generated, virtually-abelian subgroup of a compact Lie group.  Then the action $\Gamma \acts G$ by translations is residually finite.} \label{lie groups lemma}

\begin{proof}

We can represent $G$ as a subgroup of $U_n(\C)$ for some $n$.  Since the tangent spaces of $U_n(\C)$ can be identified with smoothly-varying subspaces of $M_n(\C)$, we can endow $U_n(\C)$ with a Riemannian metric coming from the real part of the Hilbert-Schmidt inner product, and this Riemannian metric is invariant with respect to multiplication by elements of $U_n(\C)$.  

Identify a set $F$ of free-abelian generators for a finite index subgroup of $\Gamma$ isomorphic to $\Z^k$.  Since these elements commute with each other, we can conjugate by an appropriate unitary so that the matrices representing them are diagonal and so the closure $\overline{\Z^k}<G<U_n(\C)$ is a torus.  Restricting the aforementioned Riemannian metric gives a Riemannian metric on this torus which is invariant and therefore lifts to $\R^k$.  So this torus is constructed from a lattice like the tori considered in the previous lemmas.  

Thus, if $x\in G$, the closure of $\Gamma x$ in $G<U_n(\C)$ is a finite disjoint union of isometric copies of that same torus (one for each of the finitely-many cosets).  Then as $\Gamma$ acts on $\overline{\Gamma x}<G$ by isometries, this action can be understood as permutations of these tori where each is transported isometrically into its image.  So if we take a copy of $E^n$ as in the previous lemma inside each of these tori, we can create an approximating action for $\Gamma\acts \overline{\Gamma x}$, as any element of $\Gamma$ will simply permute the copies of $E^n$ while applying something in $F\ltimes T^n$ to each of them.  Since our choice of $x$ at the beginning of this paragraph was arbitrary, we can repeat this process for different $x$ as necessary to create an approximating action for $\Gamma\acts G$.  \end{proof}



\prop{Suppose $\Gamma<G$ is a finitely-generated, amenable subgroup of a compact group.  Then the action $\Gamma\acts G$ by translations is residually finite.  } \label{translation actions prop}

\begin{proof}
The representation theory of compact groups shows that every compact group is an inverse limit of compact Lie groups (see for instance \cite[Chapter 16]{dixmier}).  We can assume further that the connecting maps are surjective: if they are not, replace the original limit $\displaystyle{\lim_{\leftarrow i}G_i}$ with $\displaystyle{\lim_{\leftarrow j} (\prod_{j\leq i} G_j)}$, taking the connecting maps to be projections.  Moreover, every finitely-generated, amenable subgroup of a compact Lie group is virtually abelian.  This comes from combining the Tits alternative, Lie's theorem, and Engel's theorem.  See, for instance, \cite{stack}.  

So we have that the action $\Gamma\acts G$ is conjugate to $\Gamma\acts \displaystyle{\lim_{\leftarrow} G_i}$ where $G_i$ are compact Lie groups and the latter action comes from an action of $\Gamma$ on each $G_i$ via a map $\Gamma\to G_i$ so in particular the action on each $G_i$ is implemented by an amenable (hence virtually abelian) subgroup of a compact Lie group.  By the previous lemma, each action $\Gamma\acts G_i$ is therefore residually finite.  We will show this implies residual finiteness of $\Gamma\acts \displaystyle{\lim_{\leftarrow} G_i}$, which will finish the proof.  

The pullbacks of open sets in the $G_i$ form a basis for the topology on $G$, so for any neighborhood of the diagonal in $G\times G$, we can find a neighborhood contained in it which is a finite union of products of pullbacks of balls in various $G_i$'s.  Then we can find a $K$ and $\delta>0$ such that if $(g_i)$ and $(h_i)$ are nets in the limit with $d(g_i, h_i)<\delta$ for $i\leq K$, then $((g_i), (h_i))$ is in this neighborhood of the diagonal in $G\times G$.  

Now find a $(\delta, F)$-approximating action $\Gamma\acts E\subset G_K$ and embed $E$ inside $\lim_{\leftarrow} G_i$ by taking any choice of pullbacks.  This gives a residue for the given neighborhood of the diagonal.  \end{proof}




\normalfont
%

\remark{Actions by free groups which preserve a measure of full support were previously shown in \cite[5.1]{kerr2011residually} to be residually finite, so in particular this covers the case of isometric actions by free groups.  It seems plausible therefore that the amenability assumption in the proposition above could be dropped.  }
\label{remark}

\normalfont

\normalfont
We can now extend this result to actions on compact spaces.  A special thanks to Rufus Willett for pointing out the following argument.  

\theorem{Suppose $\Gamma\acts X$ is a faithful, equicontinuous action by a finitely generated, amenable group on a compact space $X$.  Then $\Gamma\acts X$ is residually finite.  }

\begin{proof}
As before, we can assume $\Gamma\acts X$ is isometric because of \ref{equicontinuous remark}.  

Fix a finite generating set $F$.  Since the action is isometric, we can identify $\Gamma$ with a subgroup of the compact group $\text{Isom}(X)$.  The action $\Gamma\acts G := \text{Isom}(X)$ by translations is then residually finite.  Consider, for the uniform metric on $\text{Isom}(X)$, an $(\epsilon, F)$-residue $\Gamma\acts E$ where $E$ contains the identity.  

It suffices to verify the $\epsilon$-density condition locally, so pick any $x\in X$ and consider the orbit $Gx$ and the orbit map $\zeta: G\to Gx$.  Then for any $e_0\in E$ and $f\in F$, we have $f\cdot e_0 = e$ for some $e\in E$ with $d(e, f\cdot e_0)<\epsilon$, so $d(\zeta(f\cdot e_0), f\cdot \zeta(e_0)) = d(ee_0\cdot x, fe_0\cdot x)<\epsilon$, meaning $(Ex, \zeta)$ is an $(\epsilon, F)$-residue for $\Gamma\acts X$.   \end{proof}


\normalfont
The following corollary is now a direct application of \cite[3.4 and 3.6]{kerr2011residually}.  

\corollary{Suppose $\Gamma\acts X$ is a faithful, equicontinuous action by a finitely generated amenable group on a compact space $X$.  Then $C(X)\rtimes \Gamma$ is quasidiagonal.  } \qed

\remark{It should be noted that this result on quasidiagonality could already be shown, even in the case where $C(X)$ is not separable, as a consequence of \cite[Corollary B]{quasidiagonality} (together with some other results).  The proofs here are more direct, however.  Moreover, if \ref{translation actions prop} can be proved without the amenability assumption, then so can the theorem above.  It seems at least plausible therefore that all faithful isometric actions could be residually finite, which would imply results about MFness of crossed products.  }

\break
\normalfont
\bibliography{MyBibliography2-2}
\bibliographystyle{plain}

\end{document}

%% file: preamble.tex
\usepackage{amsmath}
\usepackage{amsthm}
\usepackage{mathrsfs}	
\usepackage{amssymb}
\usepackage{amscd}
\usepackage{setspace}	
\usepackage{tikz}		
\usepackage{tikz-cd}
\usetikzlibrary{matrix}	
\usepackage{graphicx}	
\graphicspath{ {pics/} }	
\usepackage{wrapfig}	
\usepackage{enumerate}	
\usepackage{bbm}		

\usepackage[utf8]{inputenc}	
\usepackage[english]{babel}	
\usepackage{blindtext}		
\usepackage[affil-it]{authblk}	

\usepackage{faktor}			
\usepackage{xparse}			

\usepackage{hyperref}		
\usepackage{cite}			

\DeclareDocumentCommand{\rfaktor}{m O{-0.5} m O{0.5}}{
  	\raisebox{#2\height}{\ensuremath{#1}}
 	 \mkern-5mu\diagdown\mkern-4mu
 	 \raisebox{#4\height}{\ensuremath{#3}}
	}

\DeclareMathOperator{\Aut}{Aut}

\newcommand{\Z}{\mathbb{Z}}
\newcommand{\R}{\mathbb{R}}

\newcommand{\C}{\mathbb{C}}

\newcommand{\acts}{\curvearrowright}

\newtheorem{theorem}{Theorem}[section]     
\newtheorem{corollary}[theorem]{Corollary}
\newtheorem{lemma}[theorem]{Lemma}

\newtheorem{prop}[theorem]{Proposition}

\newtheorem{remark}[theorem]{Remark}

\theoremstyle{definition}
\newtheorem{definition}[theorem]{Definition}

\newcommand{\Beginproof}{\begin{proof}}

\newcommand{\Endproof}{\end{proof} \normalfont}